\documentclass[]{article}

\usepackage{graphicx} % for pdf, bitmapped graphics files
\usepackage[a4paper,left=35mm,right=35mm,top=30mm,bottom=30mm,marginpar=25mm]{geometry}
\usepackage{amsmath}
\usepackage{amsthm}
\usepackage{amssymb}
\usepackage{xcolor} %[monochrome] 
\usepackage[displaymath,mathlines]{lineno}

\allowdisplaybreaks

\newcommand{\R}{{\mathbb{R}}}

\DeclareMathOperator{\T}{\mathcal{T}}

\DeclareMathOperator{\N}{\mathbb{N}}

\DeclareMathOperator*{\argmin}{arg\,min}

\newcommand{\I}{{\mathcal I}}

\newcommand{\cG}{{\mathcal G}}

\newcommand{\cI}{{\mathcal I}}

\newcommand{\cN}{{\mathcal N}}
\newcommand{\cS}{{\mathcal S}}

\newcommand{\cU}{{\mathcal U}}

\newcommand{\cT}{{\mathcal T}}

\newcommand{\cQ}{{\mathcal Q}}

\renewcommand{\phi}{\varphi}

\numberwithin{theorem}{section}
\author{Adriano Festa \thanks{Institut National des sciences appliqu\'ees (INSA) Rouen, Laboratoire de Math\'ematiques, 685 Avenue de l'Universit\'e, 76800 Saint-\'Etienne-du-Rouvray. {\tt\small adriano.festa@insa-rouen.fr} } \and Simone G\"ottlich \thanks{University of Mannheim, Department   of   Mathematics, 68131   Mannheim, Germany 
                {\tt\small goettlich@uni-mannheim.de} }
}

\title{\Large \bf A Mean Field Games approach for multi-lane traffic management \thanks{This work was partially supported by the Haute-Normandie Regional Council via the M2NUM project and the DFG project GO 1920/4-1.}}

\begin{document}

\maketitle
\date{}
% % % % % % % % % % % % % % % % % % % % % % % % % % % % % % % % %

\begin{abstract}               
In this work, we discuss a Mean Field Game approach to traffic management on multi-lane roads. The control is related to the optimal choice to change lane to reach a desired configuration of the system. Such approach is particularly indicated to model self-driven vehicles with complete information on the state of the system. The mathematical interest of the problem is that the system of partial differential equations obtained is not in the classic form, but it consists of some continuity equations (one for each lane) and a variational inequality, coming from the Hamilton-Jacobi theory of the hybrid control. We propose a consistent semi-Lagrangian scheme for the approximation of the system and we discuss how to improve its efficiency with the use of a policy iteration technique. We finally present a numerical test which shows the potential of our approach.
\end{abstract}

%===============================================================================

% % % % % % % % % % % % % % % % % % % % % % % % % % % % % % % % % % % % % % % % % % %

% % % % % % % % % % % % % % % % % % % % % % % % % % % % % % % % %

% % % % % % % % % % % % % % % % % % % % % % % % % % % % % 

\section{Introduction}
Understanding the mechanics of traffic flow systems, modeling, and efficiently predicting the future states of them is a big challenge in the applied mathematics and civil engineering. The repercussions of this theoretic issue on real-world applications are countless and socially influential. We mention the traffic forecast and the development of intelligent systems for traffic management that typically have the purpose of avoiding congestion and related costs, improving the efficiency of the road network in term of traffic fluidity, and stabilizing the traffic flow system.

Since the 1950s, many mathematical models for vehicular traffic have been proposed and studied. The largest part of them can be classified into two main categories: \emph{microscopic models} -- where every vehicle is modeled and described as a \emph{particle} -- and \emph{macroscopic} ones -- where the \emph{density} of the vehicles is taken into account.
We refer to the survey paper \cite{bellomo:2011} and to the recent monograph \cite{garavello2016models} for a general discussion about the models available in the literature. In the most classic macroscopic models, the area to which this paper is addressed, the evolution of the traffic flow is described by non-linear hyperbolic conservation laws. This choice is natural if we consider that the total number of vehicles does not change during the evolution of the system.

\begin{figure}[!t]
\begin{center}
\includegraphics[width=8cm]{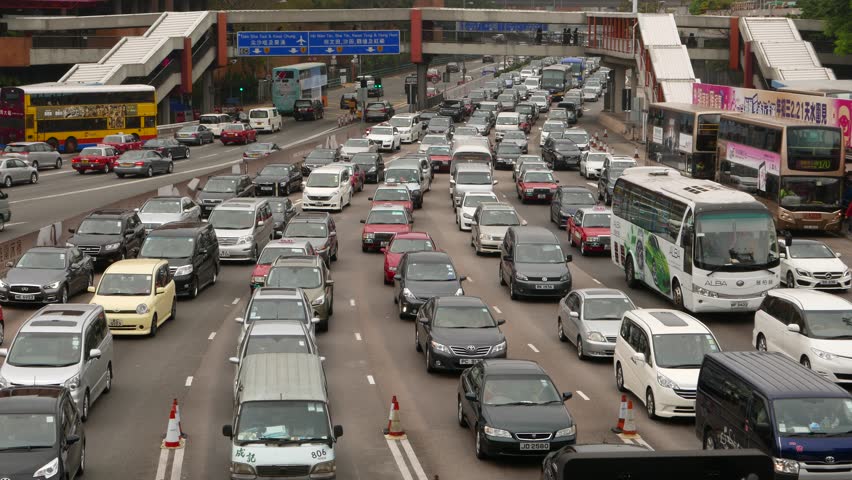}\\
\vspace{0.5cm}
\includegraphics[width=8cm]{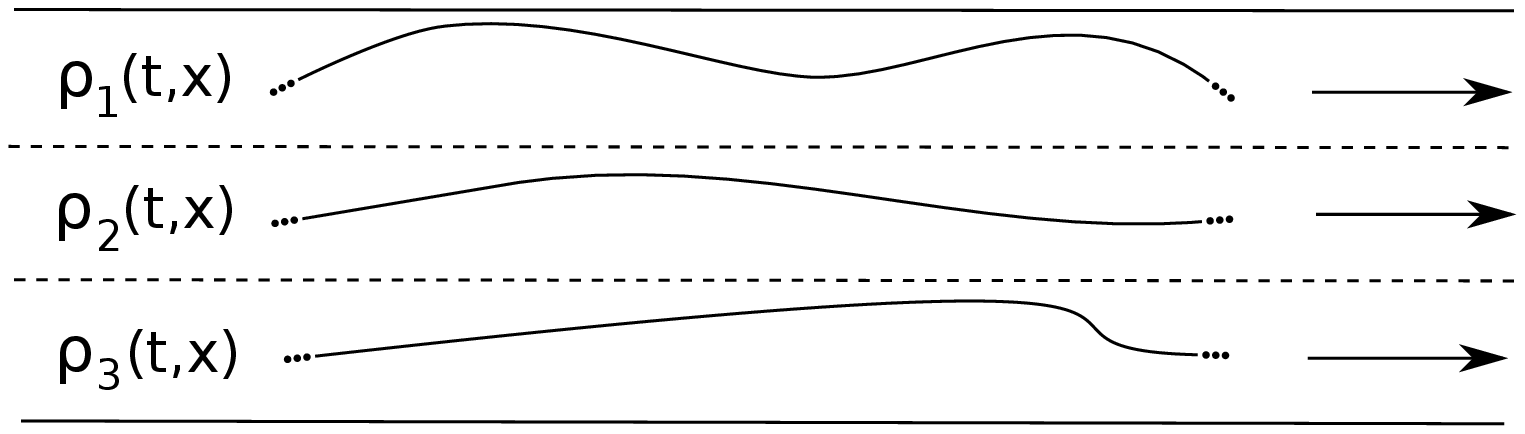}
\caption{A real world picture of multi-lane traffic and a scheme of a macroscopic model.}
\label{F:intro}
\end{center}
\end{figure}

In recent years, the availability of real-time data and the development of driver assistance and automatic driving systems intensified the request for traffic models able to perform some optimal decision making. In \cite{ll3} and \cite{Caines2}, the authors have shown that the study of strategic decision making in very large populations of small interacting individuals can be modeled by coupling a conservation law in the form of a continuity equation with a Hamilton-Jacobi equation. Nowadays, these systems are called Mean Field Games (MFG). The proposal encouraged important development in many areas, including industrial economics, algorithmics on graphs and networks, non-linear statistics, modeling of commodities price dynamics, order book dynamics, market microstructure, and crowd management.

This paper aims to use the framework of MFG to model the optimal management of a macroscopic multi-lane traffic model. In such case, a density of vehicles moves along several traffic lanes with the possibility to change lane paying a cost related to such maneuver. The aim of the control that we apply to the system leads to a desired configuration. Of course, the models require complete information about the current state of the system and related forecasting capacity by the agents (cf. \cite{festa2015collision} for collision avoidance principles). This is the typical scenario if we consider the case of a large number of interconnected vehicles. We remark that some MFG principles have been already used in the traffic-related works by \cite{cristiani} and \cite{chevalier2015micro} with a different finality and different mathematical instruments.

Finally, we observe that due to the continuous/discrete structure of the problem, the Hamilton-Jacobi equation coupled with a collection of continuity equations has the structure of a variational inequality, as studied and used in some works related to hybrid control \cite{bensoussan1997hybrid,dharmatti2005hybrid,ferretti2017hybrid}. This coupling between a series of continuity equations and the hybrid control related Hamilton-Jacobi equation is, to our knowledge, new in the literature and its analytic and numerical study is engaging and promising. We underline that the objective of this short paper is not to produce an exhaustive analysis of the framework that will be addressed in a forthcoming longer paper, but instead to show the potential of the model and the tools used in an unconventional applied scenario. For this reason, we only focus on the formal derivation of the model, and we propose a numerical approach for the approximation of the solutions. A test section shows a numerical test, and the results found. 

\section{Multi-lane traffic models}
\begin{figure}[!t]
\begin{center}
\includegraphics[width=7.5cm]{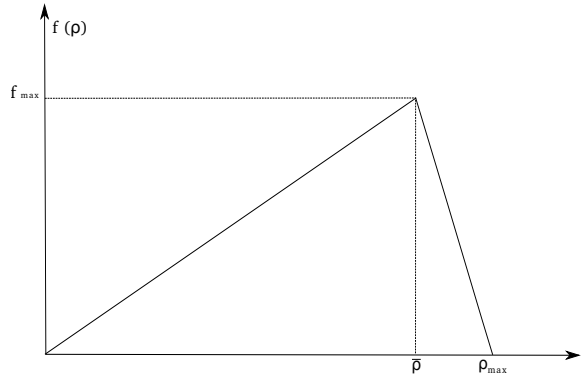}
\caption{A classic choice of piece-wise linear flux function: the maximal flux is reached for a density value $\bar\rho$. After that, the flux linearly goes to zero to prevent reaching the maximal density $\rho_{max}$.}
\label{prob}
\end{center}
\end{figure}

In order to model -- in a macroscopic way -- the traffic in a road with $n\in \N$ lanes we need to define $n$ continuity equations on $n$ domains $\R\times (0,T)$. 
 The traffic density on the $\alpha$-lane is $\rho_\alpha:\R\times (0,T)\rightarrow \R$ and follows
  \begin{equation}
\begin{cases}\label{multi}
\displaystyle  (\rho_\alpha)_t+ \left(f(\rho_\alpha)\right)_x=g(\rho_{\alpha-1},\rho_{\alpha},\rho_{\alpha+1})\\
\displaystyle \rho(x,0)=\rho_0(x)
\end{cases}
\end{equation} 
for $\alpha\in \I:=\{1,...,n\}$. Here, a source-pit term appears in the equation in order to state a link between the various lanes.\\
This is a first-order model like the LWR model introduced independently by \cite{lighthill1955kinematic} and \cite{richards1956shock}. The function $f:\R_+\rightarrow [0,f_{max}]$ gives the flux as a function of the density and it is known as \emph{fundamental diagram} in the traffic flows theory. The fundamental diagram can be modeled with a concave, piece-wise linear function as 
\begin{equation}\label{flux}
f(\rho)=\min(a\rho,b(\rho_{max}-\rho)), \quad a,b\in \R_+.
\end{equation}
In this choice, the maximal flux $f_{max}=\frac{ab}{a+b}\rho_{max}$ is reached for the density $\bar \rho=\frac{b}{a+b}\rho_{max}$ (cf. Fig. \ref{prob}).

A choice for the function $g$ proposed by \cite{shvetsov1999macroscopic} in the case of an uncontrolled system is 
\begin{multline}\label{ghelb}
g(\rho_{\alpha-1},\rho_{\alpha},\rho_{\alpha+1})=\\
\left(\frac{1}{T^L_{\alpha-1}}f(\rho_{\alpha-1})-\frac{1}{T^R_{\alpha}}f(\rho_\alpha)\right)(1-\delta_{\alpha,1})
+\left(\frac{1}{T^R_{\alpha+1}}f(\rho_{\alpha+1})-\frac{1}{T^L_{\alpha}}f(\rho_\alpha)\right)(1-\delta_{\alpha,n}).
\end{multline}
Here, the terms $T^L_\alpha$ and $T^R_\alpha$ model the tendency to change lane from lane $\alpha$ to left ($\alpha-1$) or right ($\alpha+1$). To model the boundary lanes $\alpha=1,n$  (where the change of lane can happen only in one direction) the authors use the usual Kronecker delta.

Since the exchange rates correspond to a source and a pit (or vice-versa) in consecutive lanes, the total mass of the system is preserved. This system has also been studied in \cite{klar1998hierarchy}, where the authors discussed its derivation starting from a kinetic framework. In \cite{colombo2006well} the existence and uniqueness of the solutions have been derived.

In the following, we adopt \eqref{multi} as the model for the multi-lane structure as the system to control. We point out that in the model that we propose the terms $1/T^L_\alpha$ and $1/T^R_\alpha$ disappear when we substitute them with some control functions. With this idea in mind, we need to introduce some concepts derived from hybrid systems theory.

\section{Variational inequalities and hybrid control}
We need to build an optimal control problem suitable for the multi-lane case. Since the dynamics of the system \eqref{multi} allows to switch between lanes, a natural framework to model this behavior is the one related to \emph{hybrid systems}. We summarize some concepts about the setting. We refer here to a simplified version of the one proposed in \cite{bensoussan1997hybrid} (similar formulations for the deterministic case have been proposed in \cite{branicky1998unified, dharmatti2005hybrid}). Recall that $\I=\{1,2,\ldots,n\}$ is finite, and consider the controlled system $(y,Q)$ described by:
\begin{equation} \label{eq_stato}
 \begin{cases}
 \dot y(t)=h(y(t),Q(t),u(t))\\
 y(0)= x, \ Q(0^+)=\alpha,
 \end{cases}
\end{equation}
where $x\in \R$, $\alpha\in \I$. Here, $y(t):[0,T]\rightarrow \R$ and $Q(t):[0,T]\rightarrow \I$ denote respectively the continuous and the discrete component (the lane number) of the state at time $t$. In order to end up with a Dynamic Programming equation, we also  assume that $h$ depends on $t$ only via $y$, $Q$ and $u$. The function $h :\R\times \I\times U \to \R$ represents the continuous dynamics, for a set of continuous controls given by:
$$
\mathcal{U}=\{u:(t,T) \to U \> | \> u \text{ measurable}, \ U \mbox{ compact} \},
$$
and we assume $f$  to be globally bounded and uniformly Lipschitz continuous w.r.t.~$x$.

The term $Q(t)$ takes values in the set of piece-wise constant discrete controls (or \emph{switch functions}) $\cQ$, that is:
\begin{equation}\label{switch}
\cQ = \left\{Q(\cdot):(t,T) \to \I \> |\; \> Q(t)=\sum_{i=1}^Z w_i \chi_{i}(t) \right\},
\end{equation}
where $\chi_i(t)=1$ if $t\in [t_i,t_{i+1})$ and $0$ otherwise,  $\{t_i\}_{i=1,...,Z}$ are the (ordered) times at which a switch occurs, and $\{w_i\}_{i=1,...,Z}$ are values in $\I$. 

The trajectory starts from $(x,\alpha)\in \R\times \I$. The choice of the control strategy defined as $\mathcal{S}:=\left(u,\{t_i\},\{Q(t_i^+)\}\right)$ has the objective of minimizing the following cost functional of minimum time type:
\begin{equation}\label{J}
 J(t,x,\alpha;\mathcal{S})  := \int_t^{T}\ell(y_{t,x,\alpha}(s),Q(s)) ds 
 + \sum_{i=1}^Z C\left(y_{t,x,\alpha}(t_i),Q(t_i^-),Q(t_i^+)\right),
\end{equation}
where  $\ell:\R\times \I\rightarrow \R$ is the \emph{running cost} of the trajectory $y_{t,x,\alpha}$ (solution of \eqref{eq_stato} starting at time $t$ in the point $(x,\alpha)$ and in the lane $Q(s)$ at time $s\in(t,T)$). $C:\R\times \I \times \I  \to \R_+$ is the switching cost between the dynamics, which is assumed to have a strictly positive infimum, to be bounded and Lipschitz continuous w.r.t.~$x$ and to satisfy the condition
\begin{equation}\label{triang}
C(x,\alpha_1,\alpha_2) < C(x,\alpha_1,\alpha_3) + C(x,\alpha_3,\alpha_2),
\end{equation}
for any triple of indices $\alpha_1$, $\alpha_2$ and $\alpha_3$. The condition \eqref{triang} can be seen as a strict triangular inequality for the switching cost function between the states.

The value function $V$ of the problem is then defined, for $\mathcal{S}\in\cU\times\R_+^Z\times \I^Z$, as:
\begin{eqnarray} \label{f_valore}
 V(t,x,\alpha):= \inf_\mathcal{S} J(t,x,\alpha;\mathcal{S}),
\end{eqnarray}
and is characterized via a suitable Hamilton-Jacobi-Bellman (HJB) equation. Continuity of the value function, which must be guaranteed in the classic theory of viscosity solutions is not an easy task in deterministic hybrid control problems (for a precise statement of the hypothesis we refer to \cite{dharmatti2005hybrid}), even though some more general results have been proved in a weaker framework (see \cite{bensoussan1997hybrid}). 

Through an \emph{ad hoc} adaptation of the Dynamic Programming Principle we can prove that the value function of the problem solves a HJB equation in a Quasi-Variational Inequality form. In other words, defining for $x,p\in\R$ and $\alpha\in\I$ the Hamiltonian function by
\begin{equation}\label{Ham}
 H(x,\alpha,p) := \sup_{u\in U}\{ - h(x,\alpha,u)\cdot p - \ell(x,\alpha) \}
\end{equation}
and the controlled switching operator $\cN$ by:
\begin{equation*}
 \cN \phi(x,\alpha) := \inf_{\beta\in \I} \{\phi(x,\beta)+C(x,\alpha,\beta)\},
\end{equation*}
for every function $\phi:\R\times \I$, we have a HJB equation of the following form:
  \begin{equation}\label{hjb}
 \max\left(V-\cN V,  V_t+H(x,\alpha, V_x) \right)= 0,
  \end{equation}
defined on $[0,T]\times \R\times\I$, i.e., a system of Quasi-Variational Inequalities, complemented with the boundary condition
$$
V(T,x,\alpha)=V_0^\alpha \quad x\in\R,\quad \alpha\in \I.$$

In \eqref{hjb} there are contained two separate Bellman operators which respectively provide the best possible switching and the best possible continuous control. The argument attaining the maximum in \eqref{hjb} represents the overall optimal control strategy.

\section{The mean field game system}

In this section, we use the elements described in the previous parts to build the strategic system of a mean-field type. We consider the speed of a vehicle along the lane bounded to the set $[0,f_{max}]$ so then the dynamics of the hybrid system are $h(x,\alpha,u)=uf(\rho_\alpha(x))$ and the control set $U=[0,1]$. It is clear that $h(x,\alpha,u)\geq 0$ for every choice of the parameters so then a vehicle cannot choose the direction of evolution along a road but only its speed till a maximum determined by the local density of the vehicles.

We model the switching cost between the lanes $\alpha$ and $\beta$ as the distance between them multiplied by a strictly positive parameter $\kappa\in \R_+$
\begin{equation}\label{parC}
C(x,\alpha,\beta)=\kappa \,|\alpha-\beta|.
\end{equation}
The positivity of $\kappa$ plays the important role of avoiding a chattering behavior the lanes, and it is necessary to guarantee the existence of a switching function $\cQ$.
We notice that due to the feedback nature of such switch function (implicitly stated when we used a dynamical programming principle to get a HJB), the switch control $Q$ depends on the current state of the system, then 
\begin{equation}
Q(t,x,\alpha)=\argmin_{\beta\in\cI}\left(V(t,x,\beta)+C(x,\alpha,\beta)\right).
\end{equation}
We choose the running cost function $\ell$ in order to penalize the regions of high vehicles density. This choice mirrors the one done in \cite{h02} in the different context of crowd motion models
$$\ell(x,\alpha)=\frac{1}{\max(\rho_{max}-\rho_\alpha(x),\varepsilon)},$$
for a small positive $\varepsilon$.

We are ready now to describe the complete MFG system that is
  \begin{equation}\label{MFGmain}
\begin{cases}
\displaystyle  \left(\rho_\alpha(t,x)\right)_t-\left(V_x(t,x,\alpha)f(\rho_\alpha(t,x))\rho_\alpha(t,x)\right)_x\\
\displaystyle \hspace{4.1cm}=-\chi_{Q(t,x,\alpha)\neq\alpha}f(\rho_{\alpha}(t,x))+\sum_{\beta\in \cI\setminus\{\alpha\}}\chi_{Q(t,x,\beta)=\alpha}f(\rho_{\beta}(t,x)),\\
\displaystyle \max\left[V-\cN V,\;V_t(t,x,\alpha)+ \sup_{u\in U}\left\{-u\,f(\rho_\alpha)V_x(t,x,\alpha)\right\}-\frac{1}{\max(\rho_{max}-\rho_\alpha(t,x,\alpha),\varepsilon) }\right]= 0,\\
Q(t,x,\alpha)=\argmin_{\beta\in\cI}\left(V(t,x,\beta)+C(x,\alpha,\beta)\right)
\end{cases}
\end{equation} 
defined for $(t,x,\alpha)\in  (0,T)\times \R\times\I$ and provided with the boundary condition
 \begin{equation}
\begin{cases}
\displaystyle  \rho_\alpha(0,x)=\rho^\alpha_0(x),\\
\displaystyle V(T,x,\alpha)=V^\alpha_T(x).
\end{cases}
\end{equation} 
The forward-backward structure  of these equations is a natural characteristic of the MFG systems. The densities $\rho_\alpha$ follow the optimal strategy $\cS$ and move accordingly to the respective continuity equations (with the possibility to switch among them). At the same time, such strategy must be obtained backward in time starting from a desired state of the system that we want to reach. For this reason, a natural choice (cf. \cite{lachapelle2011mean}) of the final condition $V^\alpha_T(x)$ is the distance from a desired area of the domain $\Gamma\subset\R$: in other words 
\begin{equation}\label{eik}
V^\alpha_T(x)=\inf_{y\in\Gamma}|x-y|, \quad \hbox{for every }\alpha\in \cI.
\end{equation}

%is the solution of the eikonal equation
%\begin{equation}\label{eik}
%\left\{ 
%\begin{array}{ll}
%|v_x(x)|=1, & \quad x\in \R\setminus \Gamma,\\
%v(x)=0, & \quad x\in \Gamma.
%\end{array}\right.
%\end{equation}

\section{Numerical approximation}

The first equation in \eqref{MFGmain} is a non-linear continuity equation. We can observe that the non-linear term $V_x$ in \eqref{MFGmain} -- giving the velocity vector field of evolution of the density $\rho$ --  depends {\it non-locally} on  $\rho$.  Considering these difficulties in terms of numerical approximation schemes we propose a semi-Lagrangian (SL) scheme designed to solve this case accurately. The advantages of the choice of such scheme are: the good stability of the scheme without restrictions on the choice of the discretization parameters and the monotone property implicitly stated in the formulation of the scheme.  The scheme is an adaptation of the one proposed in \cite{CS12,CS15,carlini2016PREPRINT}. We also refer the reader to \cite{FestaTosinWolfram}, where a similar scheme has been applied  to a non-linear continuity equation modeling  a kinetic pedestrian model.

\emph{Continuity equations.}
We recall here the scheme for the case of a one dimensional non-linear continuity equation on a domain $\Omega$ with some source-pit terms:
 \begin{equation}\label{FP}
 \left\{
\begin{array}{ll}
\rho_t + ( \upsilon[\rho](t,x)\, \rho )_x = g(t,x,\rho)&\hspace{0.3cm} \mbox{in }\Omega \times (0,T), \\[6pt]                         
 \rho(\cdot,0)= \rho_0(\cdot)  &\hspace{0.3cm} \mbox{in $\Omega$}.
\end{array} \right.
\end{equation}
Here, 
$\upsilon[\rho]: \Omega \times[0,T]\to\R$ is a given smooth vector field, in our case $\upsilon[\rho]:=V_x f(\rho)$ depending on $\rho$, $\rho_0$ a smooth initial datum defined on $\Omega$.  Formally, at time $t\in [0,T]$ the solution of  \eqref{FP} is given, implicitly, by the image of the measure $\rho_0 d x $ induced by the flow $x\in \Omega \mapsto \Phi(x,0,t)$, where, given $0\leq s \leq t \leq T$, $\Phi(x,s,t)$ denotes the solution of 
\begin{equation}\label{caracteristicas_continuas}
\left\{\begin{array}{rcl}
\dot \Phi(r)&=& \upsilon[\rho](\Phi,r) \hspace{0.2cm} \forall \; r\in [s,T],\\[4pt]
\Phi(s)&=& x,
\end{array}\right.
\end{equation}
at time $t$. 

Given $M\in \mathbb{N}$, we construct a grid on $\Omega$ defined by  a set of points $\mathcal{G}_{\Delta x}=\{x_i\in \Omega,\; i=1,...,M\}$ and by a set $\T$
of intervals, whose boundaries  belong to  $\mathcal{G}_{\Delta x}$ and their maximum distance is $\Delta x>0$, which form a non-overlapping coverage of $\Omega$.\\
Given $N\in \mathbb{N}$, we define the time-step $\Delta t=T/N$ and consider a uniform  partition of $[0,T]$ given by $\cT_{\Delta t}:=\{k\Delta t,\quad k=0,\hdots, N-1\}$.\\
We discretize \eqref{FP} using its representation formula by means of the flow $\Phi$. This can be considered an extension of the \emph{characteristic method} for these kind of equation (see for details \cite{CDY}).
For any $\rho \in \R^M$, $j\in \{1,...,M\}$, $k=0,\hdots, N-1$, { we  define the discrete  characteristics as}
\begin{equation*}\label{car} 
\Phi_{j,k}[\rho]:= x_{j}+\Delta t \,\upsilon[\rho](x_j,k\Delta t).
 \end{equation*}
We call $\{{\gamma_{i}} \; ; \; i=1,...,M\}$ the set of base functions, such that $\gamma_i(x_j)=\delta_{i,j}$  (the Kronecker delta) and $\sum_i \gamma_i(x)=1$ for each $x\in \overline\Omega$.\\
 We approximate the solution  $\rho$ of the problem \eqref{FP}  by a sequence $\{\rho_{k}\}_{i}=:\rho_{i,k}$, where for each $k=0,\dots,N$ $\rho_k:\mathcal{G}_{\Delta x}\to \R$ and for each $i=1,\dots,M$, $\rho_{i,k}$ approximates
$$\frac{1}{|E_i|}\int_{E_i} \rho(k \Delta t,x) d x,\quad x\in[x_i,x_{i+1}),$$
where $|E_i|=|x_i-x_{i+1}|$. 
We compute the  discrete solution $\rho_{i,k}$ by the following explicit scheme:  
\begin{equation}\label{schemefp}
\left\{
\begin{array}{ll}
\rho_{i,k+1}= G(\rho_k,i,k)+g(k\Delta t,x_i,\rho_k), & \forall k=0,..., N-1, \; \; i=1,...,M, \\[8pt]
\rho_{i,0}=\frac{ \int_{E_{i}}\rho_{0}(x) d x}{|E_i|},  &\forall i=1,...,M,
\end{array} \right.
\end{equation}
where %$m_k:=(m_{i,k})_{i=1,...,M}$ and 
 $G$ is defined by
\begin{equation*}\label{definicionG}
 G (w,i,k) :=   \sum_{j=1}^M
\gamma_{i} \left(\Phi_{j,k}\left[ w \right]\right)  w_j \frac{|E_j|}{|E_i|},
\end{equation*} 
for every $w\in \R^{M}$.

\emph{Hamilton Jacobi equation.}
We build a semi-Lagrangian scheme to give an approximation of the HJB equation of \eqref{MFGmain}. 
We consider the discrete grid of nodes $(x_j,\alpha),$ where the nodes $x_j\in\cG_{\Delta x}$ as stated previously and $\alpha\in \cI$ the lane number. In what follows, we  denote the discretization steps in compact form by $(\Delta t, \Delta x)$ and the approximate value function by $V^\Delta$.

Following \cite{ferretti2015monotone}, we write the scheme at $(x_j,\alpha)\in \cG_{\Delta x}\times \cI$ in a fixed point form as
\begin{equation}\label{scheme1}
V^\Delta(x_j,\alpha) = 
%\begin{cases}
\min \left\{ \Psi (x_j,\alpha,V^\Delta), \Sigma(x_j,\alpha,V^\Delta) \right\}. %\quad  (x_j,\alpha)\in \cG_{\Delta x}\times \cI %\\
%\Sigma(x_j,q,V^\Delta) & \text{else.}
%\end{cases}
\end{equation}
With this notation, a natural definition of the discrete jump operator $\Psi$ is given by
\begin{equation}\label{scheme11}
   \Psi (x_j,\alpha,V^\Delta) := \min_{\beta\in \cI}\left\{V^\Delta(x_j,\beta) + C(x_j,\alpha,\beta)\right\}.
\end{equation}
On the other hand, a standard semi-Lagrangian discretization of the Hamiltonian (see \cite{falcone2013semi}) is given by
\begin{equation}\label{scheme2}
\Sigma\left(x_j,\alpha,V^\Delta\right)= \min_{u\in U}\left\{\Delta t \> \ell(x_j,\alpha) + \mathbb{I}\left[V^\Delta\right] (y_j, \alpha)\right\},
\end{equation}
where $y_j:= x_j+\Delta t \> h(x_j,\alpha,u).$ The scheme is extended to all $x\in\R$ and $\alpha\in\I$ using an interpolation operator $\mathbb{I}\left[V^\Delta\right] (x,\alpha)$ which approximates the value of $V^\Delta$ at $(x,\alpha)$.
Note that because the scheme is strongly consistent and monotone, the convergence can be proved with the standard procedure derived from \cite{barles1991convergence}.

\subsection{Policy iteration algorithm for resolution of the discrete system}

The coupled discrete structure composed by the schemes \eqref{schemefp} and \eqref{scheme1} can be of difficult resolution. The issue is mostly due to the temporal direction of evolution of the two systems in \eqref{MFGmain}. In fact, while the continuity equation evolves forward in time, the variational inequality goes backward in time. The latter means that the system composed by \eqref{schemefp}, \eqref{scheme1} must be solved at the same time in all the grid $\cG_{\Delta x}\times \cT_{\Delta t}$. Various tools have been proposed to overcome the difficulty. While in \cite{CDY} the authors use a Newton iteration technique, in \cite{CS15}, for example, the authors use a more expansive fixed point iteration technique. 

We adopt a \emph{policy iteration} approach which contains the efficiency of the first idea and the simplicity of the second one. The technique consists in an alternate improvement of the control variable and an evaluation of the value function. 
We do not provide details, rather refer the reader to \cite{ferretti2016semi} where the authors use a similar approach for studying a controlled hybrid system.

\section{Tests}

\begin{figure}[!t]
\begin{center}
\includegraphics[width=8.8cm]{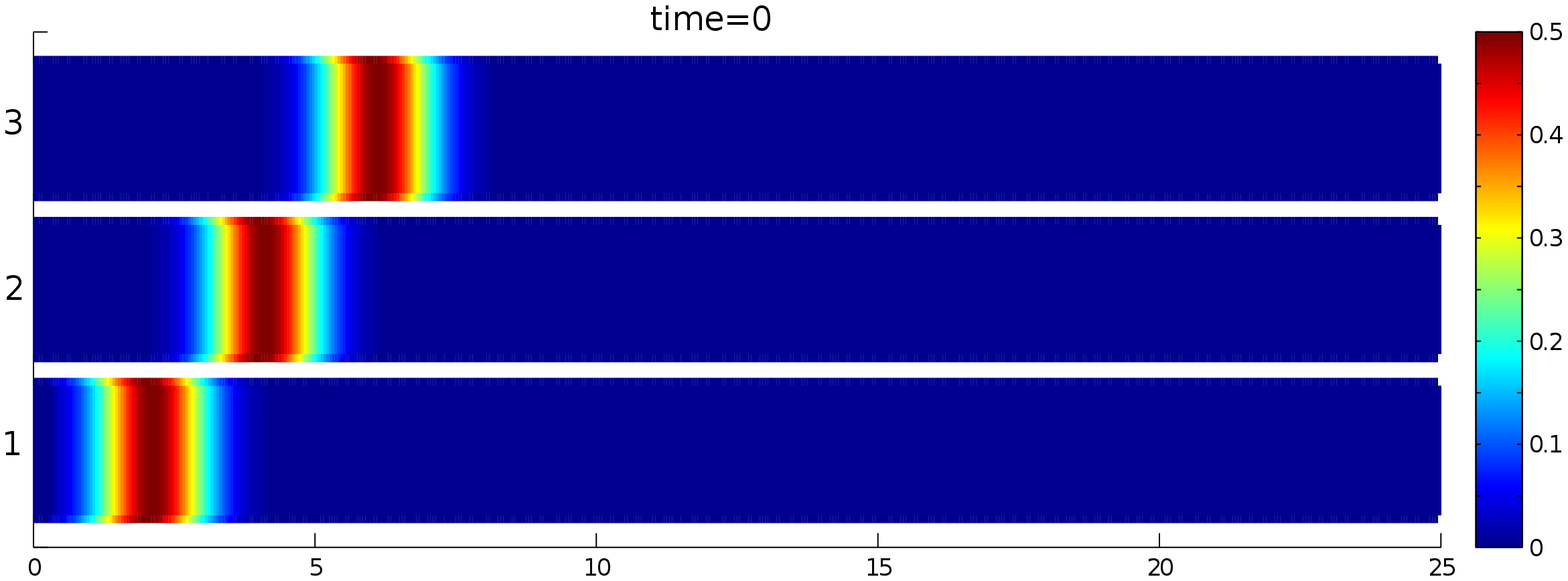}
\includegraphics[width=8.8cm]{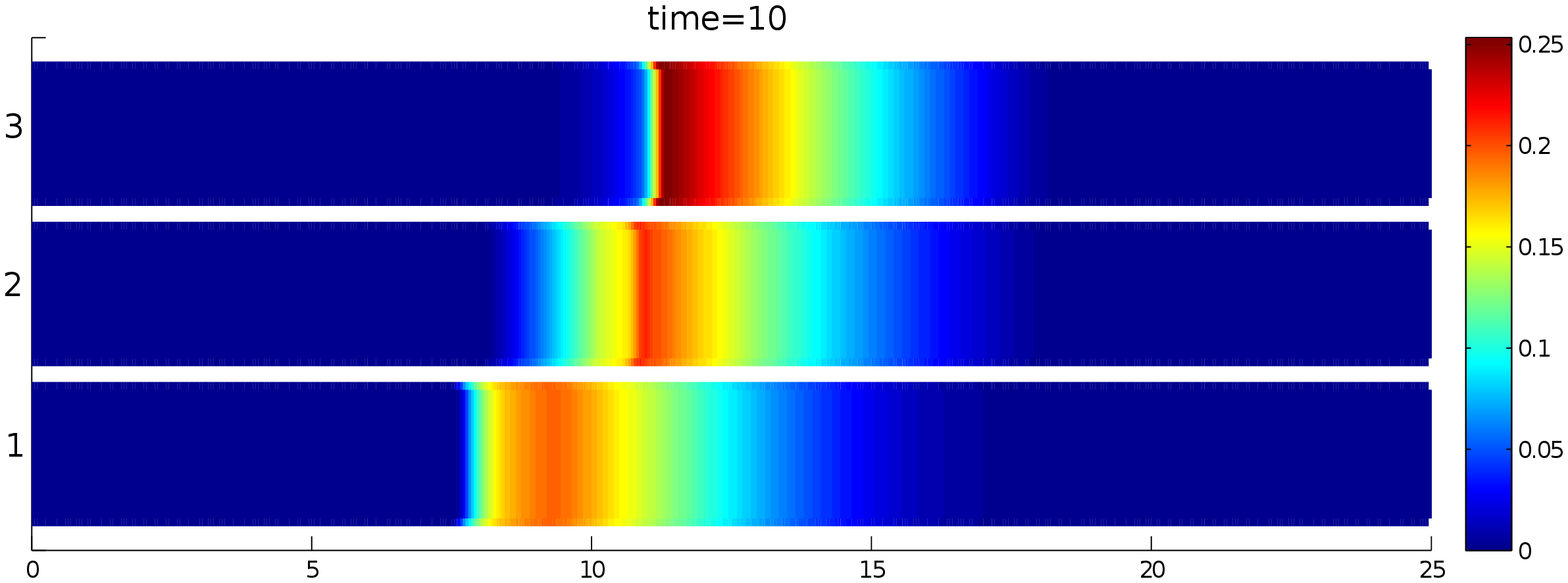}
\includegraphics[width=8.8cm]{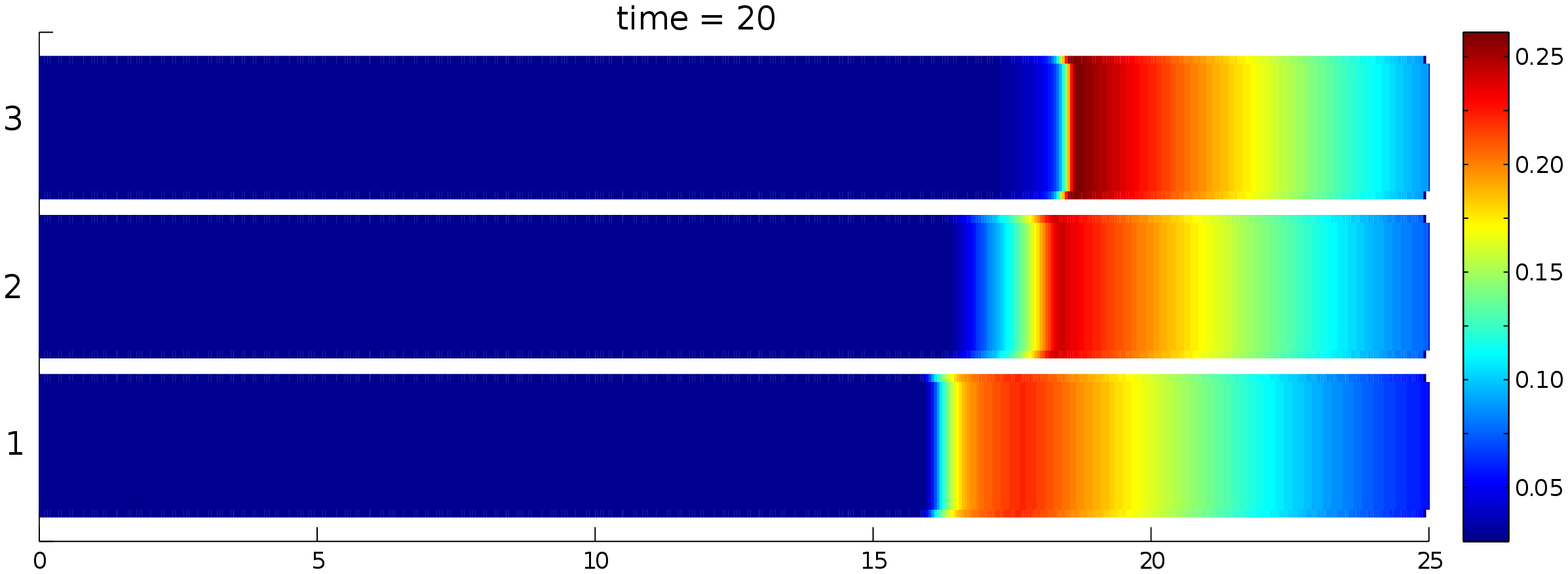}
\includegraphics[width=8.8cm]{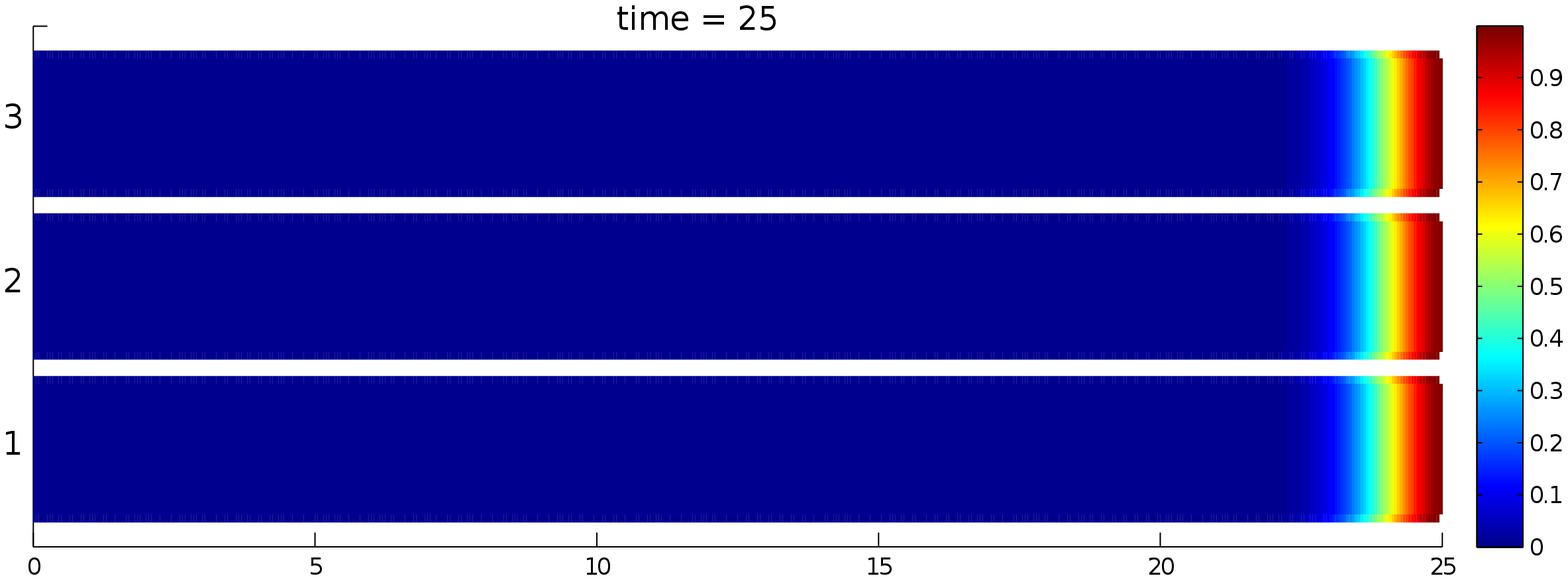}
\caption{Density distribution on the 3 lanes at various instants of the evolution of the system. The first one corresponds to the initial solution $\rho_0$, the last one is the solution at the end of its evolution at $t=25$. We underline that the space domain is $[0,25]\times\{1,2,3\}$. }
\label{F:init}
\end{center}
\end{figure}

We discuss a simple applicative scenario, where our model provides the evolution of the configurations of the system as well as the control strategies $\cS$.

We consider a $3-$lanes system, where $\Omega=[0,25]$ and the parameters of the flux function \eqref{flux} are set as $a=3$, $b=1$ and $\rho_{max}=1$. As discretization parameters we set $\Delta x=0.005$ and $\Delta t=0.01$. We remind that it is possible to adopt $\Delta t>\Delta  x$ without stability problems thanks to the well-known good features of the semi-Lagrangian scheme (see e.g. \cite{falcone2013semi}).\\
In Figure \ref{F:init} above, it is possible to see the initial configuration of the system relative to restriction on the grid of the initial data
\begin{equation*}
\rho_0^\alpha(x)=\left\{
\begin{array}{ll}
e^{-(x-2)^2}/2,& \quad\alpha=1;\\
e^{-(x-4)^2}/2,& \quad\alpha=2;\\
e^{-(x-6)^2}/2,& \quad\alpha=3.\\
\end{array} \right.
\end{equation*}

\begin{figure}[!t]
\begin{center}
\includegraphics[width=8.8cm]{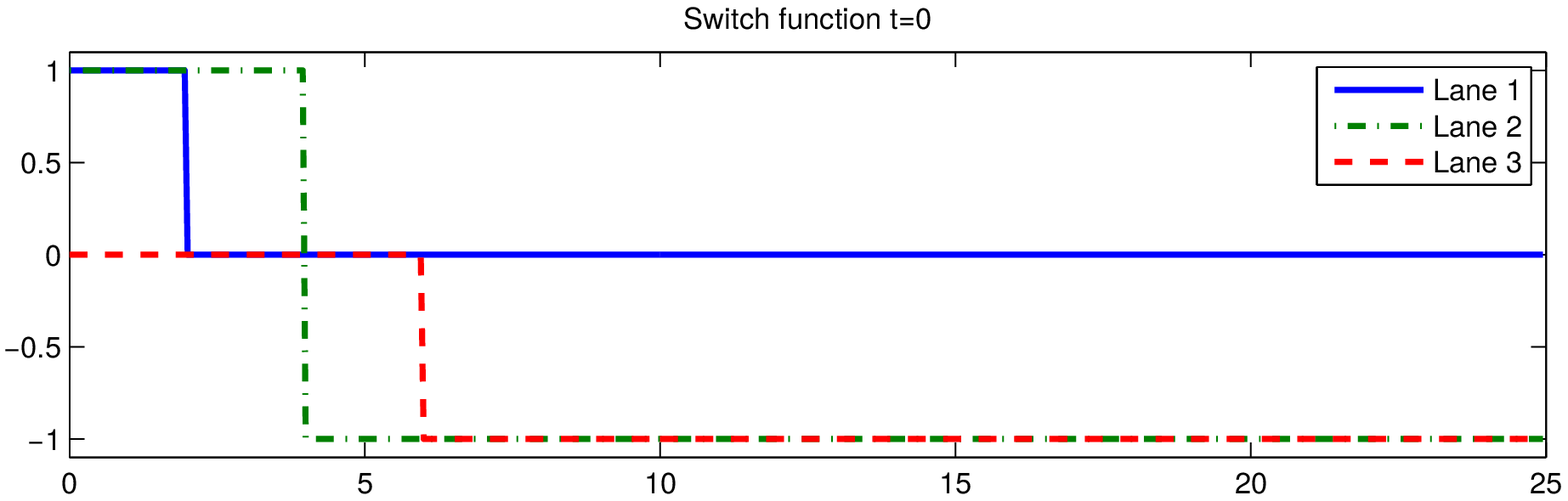}
\includegraphics[width=8.8cm]{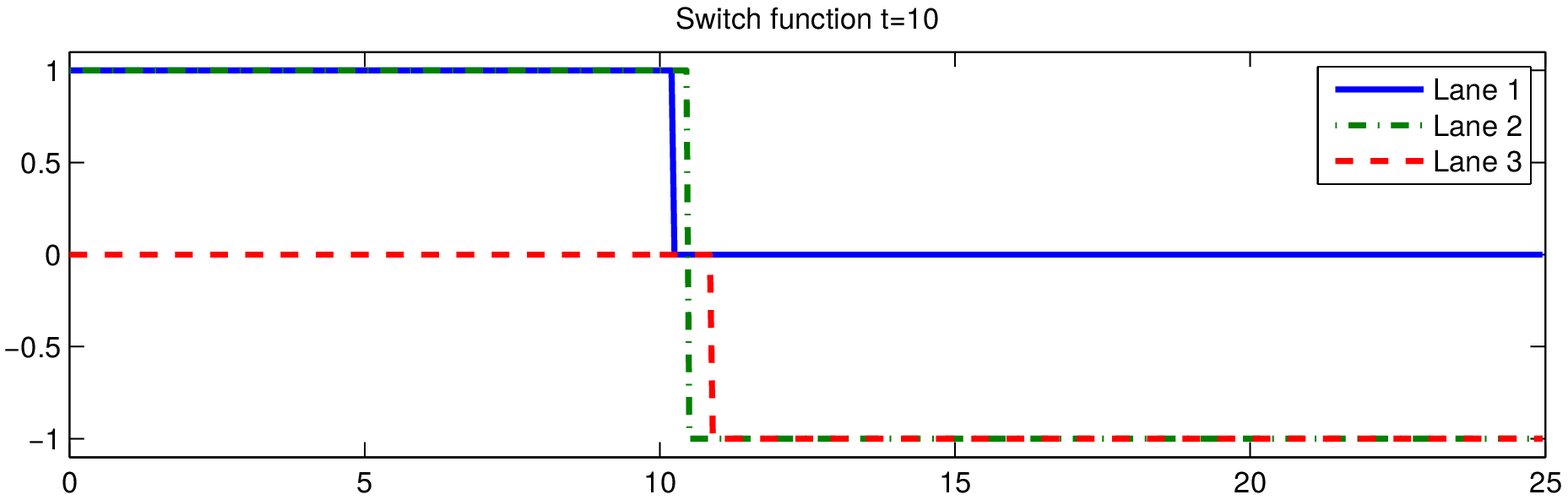}
\caption{The relative switch function $S(t,x,\alpha):=Q(t,x,\alpha)-\alpha$ in two moments of the evolution of the system. If $S(t,x,\alpha)=\pm 1$ the optimal trajectory switches respectively up/down. Above $t=0$, below $t=10$.}
\label{F:cont}
\end{center}
\end{figure}

We set the parameters of the optimization problem as $\kappa=1$ in \eqref{parC}, we substitute the control set $U=[0,1]$ with its 11 points discretization $\{0,0.1,0.2,...,1\}$ and we set $\varepsilon=10^{-5}$. We observe the behavior in the time interval $[0,25],$ where at the final time we impose some boundary conditions $V^\alpha_T(x)$ as in \eqref{eik} where $\Gamma=\{(25,\alpha), \alpha=1,2,3\}$. The initial solution has been chosen since that we are able to interpret and judge the output of the method, in term of strategies. We have mild congested areas (the maximal value of density is $0.5$) in the left part of the road moving with maximal control $u=1$ direction right. Since such densities are slightly shifted the densities tend to switch lane to maximize the occupation to reduce the congested areas and to get closer to the value $\bar \rho=1/4$ corresponding to the maximal flux (cf. Fig. \ref{F:init}, $t=10$).  This behavior can also be observed in Fig. \ref{F:cont}, where the switching functions relative to the various lanes are reported. We can observe how the lane 3, where the density is in average closer to the destination is the more congested: this is a consequence of the fact that MFG approaches look for the optimal strategy of the whole system finding a global optimum in a \emph{Nash sense}. Around the value $t=12.5$, the strategy becomes the opposite. The switching functions invert their sign, and the densities tend to redistribute uniformly on the three lanes to fulfill the final state (shown in Fig. \ref{F:init} $t=25$), where the densities, arriving at the same time on the right boundary of the domain, tend to concentrate. Consequently, they find a common configuration quite close to the maximal admissible value $\rho_{max}=1$.

\section{Conclusion}
In this paper, we have considered an MFG approach for the management of a multi-lane traffic system. The paper intends to be an initial contribution that must be fully developed in a forthcoming paper of longer extension. \\
Many questions arisen remain open at the moment. First of all the well-posedness of a system like \eqref{MFGmain} and in general of an MFG-like coupling between a system of continuity equations and a variational equation is still not proved in literature. Most importantly, is the continuity equation the adjoint of the HJB equation? In the case of a positive (as we expect) answer, and in consideration of the monotone nature of the continuity equations in \eqref{MFGmain}, a fixed point argument between the two equations may be successfully used to show existence and uniqueness of the equilibrium solution. \\
In an applied context the implementation of more complicated tests (for example with different flux parameters \eqref{flux} for different lanes) and the extension of the model to junction problems and intersecting roads is a possible development of large interest. 

\section*{Acknowledgment}
AF thanks R. Ferretti, E. Carlini, and F. Silva for the useful discussions about the MFG and the hybrid control subject.

\bibliographystyle{plain}

\end{document}